\documentclass[12pt]{article}

\usepackage{amsmath, epsfig, cite}
\usepackage{amssymb}
\usepackage{amsfonts}
\usepackage{latexsym}
\usepackage{graphicx}

\newtheorem{thm}{Theorem}[section]
\newtheorem{prop}[thm]{Proposition}

\newtheorem{cor}[thm]{Corollary}

\newtheorem{conj}[thm]{Conjecture}

\newcommand{\qed}{{\hfill\rule{4pt}{7pt}}}

\numberwithin{equation}{section}

\makeatletter \@addtoreset{equation}{section} \makeatother

\title {\bf A solution to a conjecture on the\\
 rainbow connection number\footnote{Supported by NSFC. }}

\author{
  {\small Xiaolin Chern, Xueliang Li}\\
{\small Center for Combinatorics and LPMC-TJKLC}\\
  {\small Nankai University, Tianjin 300071, P.R. China}\\
  {\small E-mail: xiaolin\_chern@yahoo.cn; lxl@nankai.edu.cn}
   }
\date{}
\begin{document}

\maketitle
\begin{abstract}
For a graph $G$, Chartrand et al. defined the rainbow connection
number $rc(G)$ and the strong rainbow connection number $src(G)$ in
``G. Charand, G.L. John, K.A. Mckeon, P. Zhang, Rainbow connection
in graphs, Mathematica Bohemica, 133(1)(2008) 85-98''. They raised
the following conjecture: for two given positive integers $a$ and
$b$, there exists a connected graph $G$ such that $rc(G)=a$ and
$src(G)=b$ if and only if $a=b\in\{1,2\}$ or $ 3\leq a\leq b$''. In
this short note, we will show that the conjecture is true.\\[3mm]
{\bf Keywords:} edge-colored graph, (strong) rainbow coloring,
(strong) rainbow connection number.\\[3mm]
{\bf AMS Subject Classification 2000:} 05C15, 05C40

\end{abstract}

\section{Introduction}

All graphs in this paper are finite, undirected, simple and
connected. We follow the notation and terminology of [1]. Let $c$ be
a coloring of the edges of a graph $G$, i.e., $c: \
E(G)\longrightarrow \{1,2,\cdots,k\},\,k\in \mathbb{N}$. A path is
called a rainbow path if no two edges of the path have the same
color. The graph $G$ is called rainbow connected (with respect to
$c$) if for every two vertices of $G$, there exists a rainbow path
connecting them in $G$. If by coloring $c$ the graph $G$ is rainbow
connected, then the coloring $c$ is called a rainbow coloring of
$G$. If $k$ colors are used in $c$, then $c$ is a rainbow
$k$-coloring of $G$. The minimum number $k$ for which there exists a
rainbow $k$-coloring of $G$, is called the rainbow connection number
of $G$, denoted by $rc(G)$.

Let $c$ is a rainbow coloring of a graph $G$. If for every pair $u$
and $v$ of distinct vertices of the graph $G$, the graph $G$
contains a rainbow $u$-$v$ geodesic (a shortest path in $G$ between
$v$ and $u$), then $G$ is called strongly rainbow connected. In this
case, the coloring $c$ is called a strong rainbow coloring of $G$.
If $k$ colors are used, then $c$ is a strong rainbow $k$-coloring of
$G$. The minimum number $k$ satisfying that $G$ is strongly rainbow
connected, i.e., the minimum number $k$ for which there exists a
strong rainbow $k$-coloring of $G$, is called the strong rainbow
connection number of $G$, denoted by $src(G)$. Thus for every
connected graph $G$, $rc(G)\leq src(G)$. Recall that the diameter of
$G$ is defined as the largest distance between two vertices of $G$,
denoted $diam(G)$. Then $diam(G)\leq rc(G)\leq src(G)$. The
following results were obtained in [2] by Chartrand et al.

\noindent \begin{prop}\label{prop1} Let $G$ be a nontrivial connected graph of
size $m$. Then\\
1. $rc(G)=1$ if and only if $src(G)=1$ .\\
2. $rc(G)=2$ if and only if $src(G)=2$.\\
3. $diam(G)\leq rc(G)\leq src(G)$ for every connected graph $G$.\qed
\end{prop}
Chartrand et al. also considered the problem that, given any two
integers $a$ and $b$, whether there exists a connected graph $G$
such that $rc(G)=a$ and $src(G)=b$ ?  and they got the following
result.

\noindent \begin{thm}\label{th1} Let $a$ and $b$ be positive integers with
$a\geq 4$ and $b \geq (5a-6)/3$. Then there exists
a connected graph $G$ such that $rc(G)=a$ and $src(G)=b$.\qed
\end{thm}

Then, combining Proposition \ref{prop1} and Theorem \ref{th1}, they
got the following result.

\noindent \begin{cor}\label{cor1} Let $a$ and $b$ be positive integers.
If $a=b$ or $3\leq a<b$ and $b\leq \frac{5a-6}{3}$, then
there exists a connected graph $G$ such that $rc(G)=a$ and $src(G)=b$. \qed
\end{cor}

Finally, they thought the question that whether the condition $b\leq
\frac{5a-6}{3}$ can be deleted ? and raised the following
conjecture:

\noindent \begin{conj}\label{con1} Let $a$ and $b$ be positive integers.
Then there exists a connected graph $G$ such that $rc(G)=a$ and
$src(G)=b$ if and only if  $a=b\in\{1,2\}$ or $ 3\leq a\leq b$. \qed
\end{conj}

This short note is to give a confirmative solution to this
conjecture.

\section{Proof of the conjecture}

\noindent {\bf Proof of Conjecture \ref{con1}:}  From Proposition
\ref{prop1} one can see that the condition is necessary. For the
sufficiency, when $a=b\in\{1,2\}$, from Corollory \ref{cor1} the
conjecture is true. So, we just need to consider the situation $
3\leq a\leq b$.

Let $n=3b(b-a+2)$, and let $H_{n}$ be the graph consisting of an
$n$-cycle $C_n: \ v_1, v_2, \cdots, v_n$ and another two vertices
$w$ and $v$, each of which joins to every vertex of $C_n$. Let $G$
be the graph constructed from $H_n$ of order $n+2$ and the path
$P_{a-1}: u_1, u_2,\cdots, u_{a-1}$ on $a-1$ vertices by identifying
$v$ and $u_{a-1}$.

First, we will show $rc(G)=a$. Because diam($G$)=$a$, by Proposition \ref{prop1}
 we have $rc(G)\geq a$. It remains to show $rc(G)\leq a$. Note that
$n=3b(b-a+2)\geq 18$. Define a coloring $c$ for the graph $G$ by the
following rules:
\[
c(e)=\left\{
\begin{array}{ll}
i & if \,\,e\,=\,\, u_{i}u_{i+1}\,\, for\,\, 1\leq i \leq a-2,\\
a-1 & if \,\,e\,=\,\,v_{i}v \,and\, i\,\, is\,\, odd,\\
a & if\,\, e\,=\,v_{i}v\,\, and\,\, i\, \,is\, even,\\
a & if\,\, e\,=\,v_{i}w\,\, and\,\,1\leq i \leq n\\
1 & otherwise.
\end{array}
\right.
\]
Since $c$ is a rainbow $a$-coloring of the edges of $G$, it follows
that $rc(G)\leq a$. This implies $rc(G)=a$.

Next, we will show $src(G)=b$. We first show $src(G)\leq b$, by
giving a strong rainbow $b$-coloring $c$ for the graph $G$ as
follows:
\[
c(e)=\left\{
\begin{array}{ll}
i & if \,\,e\,=\,\, u_{i}u_{i+1}\,\, for\,\, 1\leq i \leq a-2,\\
a-2+i & if\, \,e\,=\,\,v_{3b(i-1)+j}v \,\,for\,\,1\leq i \leq b-a+2\,\, and\,\, 1\leq j \leq 3b,\\
i & if\,\, e\,=\,\,v_{3(j-1)b+3(i-1)+k}w\,\, for\,\, 1\leq j \leq b-a+2\,\,and\,\,1\leq i \leq b\,\,\\
&and\,\,1\leq k\leq 3,\\
1 & if\,\, e\,=\,\,v_{3(i-1)+1}v_{3(i-1)+2}\,\,for\,\, 1\leq i\leq b(b-a+2),\\
2 & if\,\, e\,=\,\,v_{3(i-1)+2}v_{3(i-1)+3}\,\,for\,\, 1\leq i\leq b(b-a+2),\\
3 & otherwise\\
\end{array}
\right.
\]

It remains to show $src(G)\geq b$. By contradiction, suppose
$rc(G)<b$. Then there exists a strong rainbow $(b-1)$-coloring $c: \
E(G)\rightarrow \{1,2,\cdots,b-1\}$. For every $v_i \ (1\leq i\leq
n)$, $d(v_i,u_1)=a-1$, and the path $v_ivu_{a-2}\cdots u_{1}$ is the
only path of length $a-1$ connecting $v_i$ and $u_1$, and so
$v_ivu_{a-2}\cdots u_{1}$ is a rainbow path. Without loss of
generality, suppose
$c(u_2u_1)=1,\,c(u_3u_2)=2,\cdots,\,c(u_{a-1}u_{a-2})=a-2$. Then
$c(v_iv)\in\{a-1,a,\cdots,b\}$, for $1\leq i\leq n$. We first
consider the set of edges $A=\{v_iv,1\leq i\leq n\}$, and so
$|A|=n$. Thus there exist at least $\lceil\frac{n}{b-a+1} \rceil
\geq 3b+1$ edges in $A$ colored the same. Suppose there exist $m$
edges $v_{j_1}v,\cdots,v_{j_m}v,(1\leq j_1<j_2<\cdots<j_m\leq n)$
colored the same and $m\geq \lceil\frac{n}{b-a+1} \rceil \geq 3b+1$.
Second, we consider the set of edges
$B=\{v_{j_1}w,\cdots,v_{j_m}w\}$. Since $c(v_{j_i}w)\in
\{1,2,\cdots,b-1\}$, for $1\leq i\leq m$, then there exist at least
$\lceil\frac{m}{b-1} \rceil \geq \lceil\frac{3b+1}{b-1}\rceil \geq
4$ edges colored the same. Thus from $B$ we can choose $4$ edges of
the same color. Since $n\geq 18$, from the corresponding vertices on
the cycle $C_n$ of the four edges chosen above, we can get two
vertices such that their distance on the cycle $C_n$ is more than
$3$. Without loss of generality, we assume that the two vertices are
$v_{1}^{'},v_{2}^{'}$ and their distance in graph $G$ is $2$. Then
the geodesic between $v_{1}^{'}$ and $v_{2}^{'}$ in graph $G$ is
either $v_{1}^{'}wv_{2}^{'}$ or $v_{1}^{'}vv_{2}^{'}$. However,
neither $v_{1}^{'}wv_{2}^{'}$ nor $v_{1}^{'}vv_{2}^{'}$ is a rainbow
path. Thus the coloring $c$ is not a strong rainbow coloring of $G$,
a contradiction. Therefore $src(G)\leq b$ and so $src(G)=b$. The
proof is thus complete. \qed

\end{document}